\documentclass[12pt]{extarticle}
\usepackage[english]{babel}
\usepackage{newlfont}
\usepackage{amsbsy}
\usepackage{bbm}
\usepackage{color}
\usepackage[center]{caption2}
\usepackage{amssymb}
\usepackage{amsmath}
\usepackage{latexsym}
\usepackage{amsthm}
\usepackage{amsthm}
\theoremstyle{plain}
\newtheorem{thm}{Theorem}
\newtheorem{cor}[thm]{Corollary}
\newtheorem{lem}[thm]{Lemma}

\newtheorem{nota}[thm]{Notation}
\newtheorem{rem}[thm]{Remark}
\newtheorem{defin}[thm]{Definition}
\newcommand{\R}{\mathbb{R}}

\newcommand{\N}{\mathbb{N}}

\usepackage{mathtools}
\def\multiset#1#2{\ensuremath{\left(\kern-.2em\left(\genfrac{}{}{0pt}{}{#1}{#2}\right)\kern-.2em\right)}}
\usepackage[center]{caption2}

\begin{document}

\title{A generalization of Rohn's theorem on full-rank interval  matrices}
\author{Elena Rubei}
\date{}
\maketitle

{\footnotesize\em Dipartimento di Matematica e Informatica ``U. Dini'', 
viale Morgagni 67/A,
50134  Firenze, Italia }

{\footnotesize\em
E-mail address: elena.rubei@unifi.it}

\def\thefootnote{}
\footnotetext{ \hspace*{-0.36cm}
{\bf 2010 Mathematical Subject Classification:} 15A99, 15A03

{\bf Key words:} interval matrices, rank

{\bf \copyright} This is an Accepted Manuscript version of the following article, published in 
Linear and Multilinear Algebra:
E. Rubei "A generalization of Rohn's theorem on full-rank interval matrices" Linear and Multilinear Algebra, 68 (5), 931-939, (2020), DOI:10.1080/03081087.2018.1521366.
 It is deposited under the terms of the Creative Commons Attribution-NonCommercial-NoDerivatives License (http://creativecommons.org/licenses/by-nc-nd/4.0/), which permits non-commercial re-use, distribution, and reproduction in any medium, provided the original work is properly cited, and is not altered, transformed, or built upon in any way  }

\begin{abstract} 
 A general closed interval matrix is a matrix 
 whose entries are closed connected nonempty subsets of $\R$, while  
 an  interval matrix is defined to be a matrix 
 whose entries are closed bounded  nonempty intervals in $\R$.
 We say that a matrix $A$ with constant entries is contained in a general closed interval matrix $\mu$ if, for every $i,j$, we have that  $A_{i,j}
 \in \mu_{i,j}$.
 Rohn characterized full-rank square interval matrices, that is, square 
  interval matrices $\mu$ such that every constant matrix   contained in $\mu$ is nonsingular. In this paper we generalize this result to 
general closed interval matrices.
\end{abstract}

\section{Introduction}

Let $p , q \in \N-\{0\}$;
a    $p \times q$  interval matrix is a  $p \times q$  matrix 
whose entries are closed bounded  nonempty intervals in $\R$.
We say that a  $p \times q$  matrix $A$  with entries in $\R$ is contained in $\mu$ if, for every $i,j$, we have that $A_{i,j} \in \mu_{i,j}$. 
 
 On the other side, for any field 
 $K$, a partial matrix over $K$ is defined to be  a 
matrix where only some of the entries are given and they are elements of 
$K$;
a completion of a partial  matrix is a specification in $K$ of the unspecified entries. 
We say that a submatrix of a partial matrix is specified if all its entries are specified.

There are several papers both on interval matrices and
on partial matrices. 
On partial matrices, there is a wide literature about
the problem of determining the maximal and the minimal rank of the completions of a partial matrix.
 We quote, for instance, the papers  \cite{CD},  \cite{Woe} and 
\cite{CJRW}. In the last, Cohen,  Johnson,  Rodman and Woerdeman determined 
the maximal rank of the completions of a partial
matrix in terms of the ranks and the sizes of its maximal  
specified submatrices; see also \cite{CD}
for the proof.
 The problem of determining the minimal rank of   the completions of a partial matrix seems more difficult  and it has been solved only in some particular cases,
  see for example  \cite{Woe1} for the case of triangular matrices and the recent paper \cite{CP}. 
 
 Also  interval matrices have been widely studied and in particular
there are several papers studying when a $p \times q $ interval  matrix $\mu$ has full rank, that is when all the matrices contained in $\mu$ have rank equal to $\min\{p,q\}$. 
For any  $p \times q $ interval matrix $\mu = ( [m_{i,j}, M_{i,j}])_{i,j}$ with $m_{i,j} \leq M_{i,j}$,
let $C^{\mu}$,$\Delta^{\mu}$ and $|\mu|$ be 
 the  $p \times q$   matrices such that $$ C^{\mu}_{i,j}= \frac{m_{i,j}+ M_{i,j}}{2}, \hspace*{1cm} \Delta^{\mu}_{i,j}= \frac{M_{i,j}- m_{i,j}}{2}, \hspace*{1cm} 
|\mu|_{i,j} = \max\{|m_{i,j}|,| M_{i,j}|\} $$ for any $i,j$. 
For any nonzero natural number $p$, let $Y_p=\{-1,1\}^p$ and, for any $x
 \in Y_p$, denote by $T_x$ the diagonal matrix whose diagonal is $x$. Finally, for any $p \times p$  interval matrix 
 $\mu =( [m_{i,j}, M_{i,j}])_{i,j}$, with $m_{i,j} \leq M_{i,j}$ for any $i,j$, and 
  for any $x, y \in Y_p$,
 define the  matrix $C^{\mu}_{x,y}$ as follows:
  $$C^{\mu}_{x,y} = C^{\mu} - T_x \, \Delta^{\mu} \, T_y.$$
The following theorem
characterizes full-rank square interval matrices:

\begin{thm} {\bf (Rohn, \cite{Rohn})} 
Let   $\mu =( [m_{i,j}, M_{i,j}])_{i,j}$ be  a  $p \times p$  interval matrix, where
 $m_{i,j} \leq M_{i,j}$  for any $i,j$.
Then  $\mu$ is a full-rank interval matrix 
if and only if,  for each
$x,y ,x',y'\in Y_p$,  $$det(C^{\mu}_{x,y}) \,det(C^{\mu}_{x',y'})>0. $$
\end{thm}

See  \cite{Rohn} and  \cite{Rohn2} for other characterizations. Finally, as to interval matrices, we want to quote also 
the following theorem characterizing 
full-rank $p \times q$ interval matrices, see
\cite{Rohn3}, \cite{Rohn4}, \cite{Shary}:

\begin{thm}  {\bf (Rohn)} 
A $p \times q$ interval matrix $\mu$ with $p \geq q$ has full rank if and only if the system
of inequalities $$\hspace*{2cm} |C^{\mu} x| \leq \Delta^{\mu}
|x|, \hspace*{1.5cm} x \in \R^q$$ has only the trivial solution $x=0$. 
\end{thm}

  Obviously the problem of partial matrices in the real case and the one of interval matrices are connected; in fact  we can consider matrices whose entries
   are closed connected nonempty subsets of $\R$;  these matrices generalize both the interval matrices and the partial matrices.
   We call a matrix 
 whose entries are closed connected nonempty subsets of $\R$
  a ``general closed interval matrix''.
  
  In this paper we generalize Rohn's result 
  on full-rank square interval matrices
  to 
general closed interval matrices, see Theorem \ref{gen} for the precise statement.

\section{Notation and some recalls}

$\bullet$  Let $\R_{>0}$ be the set $\{x \in \R | \; x >0\}$ and
let $\R_{\geq 0}$ be the set $\{x \in \R | \; x  \geq 0\}$; we define analogously $\R_{<0}$ and $\R_{ \leq 0}$.

$\bullet $ Throughout the paper  let $p , q \in \N-\{0\}$. 

$\bullet $ Let $\Sigma_p$ be the set of the permutations on $\{1,....,p\}$. For any 
permutation $\sigma$, we denote the sign of $ \sigma$ by 
 $ \epsilon (\sigma)$.

$\bullet$ Let $M(p \times q, \R) $ denote the set of the $p \times q $  matrices
with entries in $\R$. For any $A
\in M(p \times q, \R)$, let $rk(A) $ denote the rank of $A$
and let $A^{(j)}$ be the $j$-th column of $A$.

$\bullet$ For any $p \times q$ general closed interval matrix $\mu$, any $\{i_1, \ldots, i_s\} \subset \{1, \ldots , p\}$ and any $\{j_1, \ldots, j_r\} \subset \{1, \ldots , q\}$, we denote by $\mu_{\widehat{i_1,\ldots , i_s}, \widehat{j_1, \ldots , j_r}} $ the matrix obtained from $\mu$  by deleting the rows $i_1, \ldots , i_s$ and the columns  $j_1, \ldots , j_r$.

\begin{defin} A {\bf general  interval matrix} is a matrix 
 whose entries are connected nonempty subsets of $\R$.
 
 A {\bf general closed interval matrix} is a matrix 
 whose entries are closed connected  nonempty subsets of $\R$.
 
 An {\bf interval matrix} is a matrix 
 whose entries are closed bounded nonempty intervals of $\R$.

 Let   $\mu $ be  a $p \times q $  general interval matrix.
 As we have already said, 
 given a matrix $A \in M(p \times q, \R) $, we say that  $ A \in \mu $ if   $a_{i,j} \in  \mu_{i,j} $ for any $i,j$.
 
      We define
$$mrk(\mu) = min\{rk(A)   | \; A \in \mu     \},$$
$$Mrk(\mu) = max\{rk(A) | \; A \in \mu     \}.$$
We call them respectively {\bf minimal rank} and {\bf maximal rank} of $\mu$.
Moreover, we define 
$$rkRange (\mu) =  \{ rk(A)  | \; A \in \mu \};$$ 
we call the set above the {\bf rank range} of $\mu$.

We say that the  entry $i,j$ of $\mu$ is a {\bf constant}
 if $ \mu_{i,j}$ is a subset of $\R$ given by only one element. 

We say that the  entry $i,j$ of $\mu$ is {\bf bounded}
 if $ \mu_{i,j}=[a,b]$ for some $a,b \in \R$; we say that  the  entry $i,j$ of $\mu$ is  {\bf half-bounded}
 if either  $ \mu_{i,j}=[a,+ \infty ) $ ({\bf left-bounded}) or 
  $ \mu_{i,j}=(- \infty ,a] $  ({\bf right-bounded})
 for some $a \in \R$.
\end{defin}

\begin{rem}
Let $\mu$ be an interval matrix. Observe that  $$rkRange(\mu) =
  [mrk(\mu)  , Mrk(\mu)] \cap \N.$$
\end{rem}

See \cite{Rub} for a proof.

\begin{defin}
Given a $p \times p $ interval  matrix, $\nu$, a {\bf partial generalized
diagonal} ({\bf pg-diagonal} for short) of length $k$ of $\nu$ is a $k$-uple of the kind $$
(\nu_{i_1, j_1},\dots, \nu_{i_k, j_k})$$ 
for some  $\{i_1, \dots i_k\}$ and $ \{j_1, \dots, j_k\} $ subsets of $ \{1,\dots ,p\}$.

Its {\bf complementary matrix} is defined to be the submatrix of $\mu$ given by the rows and columns whose indices are respectively in  $\{1,\ldots , p\}- \{i_1, \ldots , i_k\}$ and 
in  $\{1,\ldots , p\}- \{j_1, \ldots , j_k\}$. 

We say that a pg-diagonal is {\bf 
totally nonconstant} if  all its entries are not constant.

We define $det^c(\mu) $ as follows:
$$ det^c(\mu) = \sum_{\sigma \in \Sigma_p \; s.t. \; \mu_{1, \sigma(1)}, \dots, \mu_{p, \sigma(p)} \;  are \; constant} \epsilon (\sigma) \,
\mu_{1, \sigma(1)} \cdot \ldots  \cdot\mu_{p, \sigma(p)} $$ 
if there exists $\sigma \in   \Sigma_p $ such that $ \mu_{1, \sigma(1)}, \dots, \mu_{p, \sigma(p)}$    are  constant; we define $det^c (\mu)$ to be equal to $0$ otherwise.
\end{defin}

The following theorem and corollary were proved in \cite{Rub} for interval matrices; the same 
proof yields the results for general closed interval matrices:

\begin{thm} \label{Mrkquadrate}
Let $\mu $ be a $p \times p $ general closed interval matrix. Then $Mrk(\mu) < p$ if and only if the following conditions hold:

(1) in $\mu$ there is no totally nonconstant pg-diagonal of length $p$,

(2) the complementary matrix of every  totally nonconstant pg-diagonal of  length between $0$ and $p-1$ has $det^c $ equal to $0$ (in particular $det^c (\mu)=0$).

\end{thm}

\begin{cor} \label{Mrk}
Let $\mu $ be a general closed  interval matrix. Then $Mrk(\mu)$
is the maximum of the natural numbers $t$ such that there is a $ t \times t $ submatrix  of
 $\mu$ either with a totally nonconstant pg-diagonal of  length 
 $t$ or with a totally nonconstant
 pg-diagonal of length between $0$ and $t-1$ whose complementary matrix has $det^c  \neq 0$.
 \end{cor}

 \begin{nota} Let  $\mu$ be a general closed interval matrix.

$\bullet $ 
We denote by $\tilde{\mu}$ the matrix obtained from $\mu$ by replacing the entries $(-\infty , + \infty)$ with $0$.
 
 $\bullet$ 
 We denote by $ \overline{\mu}$ 
 the matrix obtained from $\mu$ by replacing every  entry of kind  $[a, + \infty)$, for some $a \in \R$,  with $a$ and every entry of kind $(-\infty, b]$, for some $b \in \R$, with $b$.
 
 $\bullet $ 
 We denote by $\mu_l$ 
  the matrix obtained from $\mu$ by replacing every entry of kind  $[a, b] $, for some $a ,b \in \R$ with $a \leq b$,  with $a$.
 We denote by $\mu_r$ 
  the matrix obtained from $\mu$ by replacing  every entry  of kind  $[a, b] $, for some $a ,b \in \R$ with $a \leq b$,  with $b$.

\end{nota}

 \begin{defin} \label{vertexmatrix}
Let $\mu $ be a  general closed interval matrix. 
  
We say that $\gamma $ is  a {\bf vertex matrix} of $\mu$ if  $\gamma_{i,j} \in \{m_{i,j}, M_{i,j}\}$ for any $i,j$ such that $\mu_{i,j}$ is a bounded interval and $\gamma_{i,j}= \mu_{i,j}$ otherwise.


We say that a  vertex matrix $\gamma$ of $\mu$ is {\bf of even type} if,
for every  $2 \times 2$ submatrix of $\mu$ such that all its entries are bounded intervals,
either  the number of the entries of the corresponding submatrix of $\gamma$ that are equal to the minimum of the corresponding  entries of $\mu$ is even or some of its entries are constant. 
We say that it is {\bf of odd type} if it is not of even type.  \end{defin}
  
  {\bf Example.} Let $$\mu=\begin{pmatrix}
[1,2] & [2,3] & [2,+\infty) \\
[-3,4] & [-1, 5] & [1,4]
  \end{pmatrix}.$$ Then $\gamma=\begin{pmatrix}
1 & 2 & [2,+\infty) \\
4 &  5 & 1
  \end{pmatrix}$ is a vertex matrix of $\mu $ of even type, while $\delta=\begin{pmatrix}
1 & 3 & [2 , +\infty)\\
4 &  5 & 1
  \end{pmatrix}$ is a  vertex matrix of $\mu$ of odd type. 

\bigskip

It is easy to see  that, given a $p \times p$ interval matrix $\mu$ and $x , y \in Y_p$, the matrix $C^{\mu}_{x,y}$ is a vertex matrix of $\mu $ of even type and every vertex matrix of $\mu $ of even type is equal to $C^{\mu}_{x,y}$ for some $x,y \in Y_p$.
So, by using Definition \ref{vertexmatrix}, we can restate Rohn's theorem as follows:

\begin{thm} {\bf (Rohn, \cite{Rohn})}  \label{Rohnbis}
Let   $\mu$ be  a  $p \times p$  interval matrix.
Then  $\mu$ is a full-rank interval matrix 
if and only if,  for any vertex matrices $A_1, A_2$ of even type of $\mu$, $$det(A_1) \,det(A_2)>0. $$
\end{thm}

\section{The main result}

\begin{rem} \label{polin}
Let $p(x_1,\ldots, x_n)$ be a polynomial
with coefficients in $\R$ of degree $1$ in every variable. Let $\overline{x} \in \R^n$.
Then $p(x) \geq 0 $ (respectively $p(x)>0$)  for every $x \in \overline{x} + \R^n_{\geq 0}$ if and only if 
$p(\overline{x}) \geq 0$ (respectively $p(\overline{x})>0$)
and $\frac{ \partial \,p }{\partial x_i} (x) \geq 0 $     for every $x 
 \in \overline{x} + \R^n_{\geq 0}$ and for every $i=1, \ldots, n$.

\end{rem}

\begin{lem} \label{lem1}
Let $\mu $ be a general closed interval $p \times p$ matrix. 
Then $\mu $ is full-rank if and only if
$\tilde{\mu}$ is full-rank and for every 
$i,j$ such that $ \mu_{i,j}= (-\infty, + \infty)$ 
we have that $ Mrk( \mu_{\hat{i}, \hat{j}})
 < p-1$.
\end{lem}

\begin{proof} Suppose 
 $\mu $ is full-rank. Then  obviously  
$\tilde{\mu}$ is full-rank.  Moreover, let 
$i,j$ be such that $ \mu_{i,j}= (-\infty, + \infty)$. Then the  determinant of every $A \in  \mu_{\hat{i}, \hat{j}}$ must be zero; otherwise, suppose there  exists  $A \in  \mu_{\hat{i}, \hat{j}}$ with $det(A) \neq 0$; then,
for any choice of
 $x_{i,s}$ in  $\mu_{i,s}$ for any $s \in \{1, \ldots , p\}-\{j\}$
  and $x_{r,j}$ in  $\mu_{r,j}$ for any $r \in \{1, \ldots , p\}-\{i\}$,
we can choose $x_{i,j} \in  \mu_{i,j}$ such that, if we define $X$ to be the matrix   
such that   $X_{r,s}= x_{r,s}$ for $(r,s)$ such that 
either $r=i$ or $s=j$   and
$X_{\hat{i}, \hat{j}}=A$, we have that 
the determinant of $X$  is zero, which is absurd since 
$\mu$ is full-rank.
So $ Mrk( \mu_{\hat{i}, \hat{j}})
 < p-1$ and we have proved the right-handed implication.
 
 On the other side, let $(i_1,j_1), \ldots, 
 (i_s, j_s) $ be the indices of the entries of $\mu $ equal to $ (-\infty, + \infty)$ and suppose that $\tilde{\mu}$ is full-rank and that $ Mrk( \mu_{\hat{i_k}, \hat{j_k}})< p-1$
 for every $k=1, \ldots ,s$. 
 Let  $B \in M(p \times p, \R)$ be such that  
$B \in   \mu$. For every  $k=1, \ldots, s$, let $B_k$ be the matrix obtained from $B$ by 
 replacing the  entries $(i_1,j_1), \ldots , (i_k, j_k) $  with $0$. Then 
 $det (B)= det (B_1)= \ldots =det(B_s)  $ 
 since $ Mrk( \mu_{\hat{i_k}, \hat{j_k}})< p-1$ for every $k=1, \ldots ,s$.
Obviously $B_s \in \tilde{\mu}$, so $det (B_s)$ is nonzero  since $\tilde{\mu}$ is full-rank; hence $det(B)$ is nonzero  and
we conclude.
\end{proof}

\begin{lem}  \label{lem2}
Let $\nu $ be  a general closed interval $p \times p $ matrix with only bounded or half-bounded entries. Then $\nu $ is full-rank
 if and only if, for any $\gamma $ vertex matrix of $\nu$ of even type, we have that $\gamma $  is full-rank and  $det(\overline{\gamma})$ has the same sign as $det( \overline{\nu}_l)$.  
\end{lem}

\begin{proof} $\Longrightarrow$ 
Obviously if $\nu$ is full-rank, then the determinant of all
the matrices contained in $\nu $ must have the same sign (since the determinant is a continous function and the image by a continous function of a connected subset is connected), in particular 
 $det(\overline{\gamma})$ must have the same sign as $det( \overline{\nu}_l)$.
 
 $\Longleftarrow$ 
Let $A \in M(p \times p, \R) $ be such that $A \in \nu$; define 
$\alpha $ to be the interval matrix such that 
$\alpha_{i,j} = \nu_{i,j}$ if $\nu_{i,j}$  is bounded and $\alpha_{i,j} = a_{i,j}$ if $\nu_{i,j}$  is half-bounded.

 To prove that $A$ is full-rank, obviously it is sufficient to prove that $\alpha$ is full-rank.
 By Rohn's theorem (see Theorem \ref{Rohnbis}), to prove that $\alpha$ is full-rank, it is sufficient to prove that  any two
 vertex matrices $A_1,A_2$   of even type of $\alpha$  are  full-rank and the sign of  their determinant is the same.
Obviously, for $i=1,2$, the matrix  $A_i$ is contained in a vertex matrix $\gamma_i$ of even type of $\nu$;  by assumption,  for $i=1,2$, the matrix $\gamma_i$ is full-rank, so the matrix $A_i$ is full-rank
and its determinant has the same sign
as  $det(\overline{\gamma}_i)$.
Moreover,   by assumption the sign of 
$det(\overline{\gamma}_i)$ is equal to the 
sign of $det( \overline{\nu}_l)$ for $i=1,2$;
 in particular $det(\overline{\gamma}_1)$
 and $det(\overline{\gamma}_2)$ have the same sign, so $det(A_1)$ and $det (A_2)$
 have the same sign, as we wanted to prove.
\end{proof}

\begin{nota}
Let $\rho$ be a square general closed interval matrix. We say that $DET(\rho)$ is 
greater than $0$ (respectively less than $0$, equal to $0$...) if,  for any $A \in \rho$, we have that 
$det(A)$ is greater than $0$ (respectively less than $0$, equal to $0$...). More generally,  given a function $f: \R^n \rightarrow \R$ for some $n \in \N-\{0\}$,
we say that $ f(x_1, \ldots, x_{n-1}, DET(\rho) ) \geq 0 $  (respectively $\leq 0$, $>0 $, $< 0$) if   $ f(x_1, \ldots, x_{n-1}, det(A) ) \geq 0 $  (respectively $\leq 0$, $>0 $, $<0$) for any $A \in \rho$.
 
\end{nota}

\begin{lem}  \label{lem3}
Let $\rho $ be  a general closed interval $p \times p $ matrix with only constant or half-bounded entries. Then $\rho $ is full-rank
 if and only if $det(\overline{\rho}) \neq 0$ and, 
 for any $(i_1,j_1), \ldots, 
 (i_s,j_s) \in \{1,\ldots , p \} \times   \{1,\ldots , p \} $ such that $i_1, \ldots , i_s$ are distinct and $j_1, \ldots , j_s$ are distinct and $\rho_{i_1,j_1}, \ldots , \rho_{i_s,j_s} $ are half-bounded, 
$$ (-1)^{i_1+j_1+ \widetilde{i_2} + \widetilde{j_2}+ \ldots + \widetilde{i_s}+ \widetilde{j_s}+ \chi((i_1,j_1),
\ldots, (i_s, j_s) ) }  det (\overline{\rho}) \,det (\overline{\rho}_{\widehat{i_1,\ldots, i_s}, \widehat{j_1, \ldots, j_s}} )  \geq 0, $$
where: 

$\bullet $ the determinant of a $0 \times 0 $ matrix is defined to be $1$,

$\bullet$ 
$\chi((i_1,j_1),
\ldots, (i_s, j_s) )$ is defined to be  the number of the  right-bounded intervals in 
$\rho_{i_1,j_1}, \ldots, \rho_{i_s,j_s}$,

$\bullet$ 
 for $t=2, \ldots, s$, we define 
$\widetilde{i_t}$ to be  $i_t$ minus the number of the elements among $i_1, \ldots i_{t-1}$ smaller than 
$i_t$ and  $\widetilde{j_t}$ to be  $j_t$ minus the number of the elements among $j_1, \ldots j_{t-1}$ smaller than  $j_t$.

\end{lem}

\begin{proof}
By Remark \ref{polin}, the matrix $\rho $ is full-rank if and only if   $det(\overline{\rho}) \neq 0$ and,
for any $i_1, j_1$ such that $\rho_{i_1,j_1}$ is half-bounded,
$$ (-1)^{i_1+ j_1 + \chi(i_1,j_1)} \, det (\overline{\rho} ) \,DET (\rho_{\hat{i_1},\hat{j_1}} ) \geq 0 $$ 
and, again by  Remark \ref{polin}, the last condition holds  if and only if, 
for any $i_1, j_1$ such that $\rho_{i_1,j_1}$ is half-bounded, we have that
$$ (-1)^{i_1+ j_1 + \chi(i_1,j_1)} \, det (\overline{\rho} ) \,det (\overline{\rho}_{\hat{i_1},\hat{j_1}} ) \geq 0 $$ 
and, for any  $i_2, j_2$ such that $i_1 \neq i_2$ and $j_1
\neq j_2$ and 
$\rho_{i_2,j_2}$ is half-bounded, we have that
$$ (-1)^{i_1+ j_1 +\widetilde{i_2} + \widetilde{j_2} + \chi((i_1,j_1), (i_2,j_2))} \, det (\overline{\rho} ) \,DET (\rho_{\widehat{i_1,i_2},\widehat{j_1,j_2}} ) \geq 0 $$ 
and so on.
\end{proof}

\begin{thm} \label{gen}
Let $\mu $ be a general closed interval matrix. 
Then $\mu $ is full-rank if and only if the following conditions hold:

(1) for every  $i,j$ such that  $\mu_{i,j} = (- \infty, +\infty) $ we have that in  $\mu_{\hat{i}, \hat{j}} $ 
 there are no totally nonconstant pg-diagonal of  length 
 $p-1$ and the complementary matrix of every   totally nonconstant
 pg-diagonal in  $\mu_{\hat{i}, \hat{j}} $ of  length between $0$ and $p-2$  has $det^c $ equal to $0$.

(2) the product of the determinant every  vertex matrix of  $\overline{\tilde{\mu}}$ of even type and the determinant of $\overline{\tilde{\mu}}_l$ 
is positive, 

(3)  for every vertex matrix $\gamma$ of $\tilde{\mu}$ of even type and every $(i_1,j_1),\ldots, (i_s, j_s) 
\in \{1,\ldots , p \} \times   \{1,\ldots , p \} $ such that $i_1, \ldots , i_s$ are distinct and $j_1, \ldots , j_s$ are distinct and 
$\mu_{i_1,j_1}, \ldots , \mu_{i_s,j_s} $ are half-bounded, we have:
$$ (-1)^{i_1+j_1+ \widetilde{i_2}+ \widetilde{j_2}+\ldots + \widetilde{i_s}+ \widetilde{j_s}+ \chi((i_1,j_1),
\ldots, (i_s, j_s) ) } 
det (\overline{\tilde{\mu}}_l) \,
det (\overline{\gamma}_{\widehat{i_1,\ldots, i_s}, \widehat{j_1, \ldots, j_s}} )  \geq 0 ,$$
where:

$\bullet $  the determinant of a $0 \times 0 $ matrix is defined to be $1$,

$\bullet$  the number $\chi((i_1,j_1),
\ldots, (i_s, j_s) )$ is defined to be the number of right-bounded 
intervals  in 
$\mu_{i_1,j_1}, \ldots, \mu_{i_s,j_s}$,

$\bullet $ 
 for $t=2, \ldots, s$,  we define 
$$\widetilde{i_t} := i_t - \sharp \{i_r\; \mbox{\em for} \; r=1, \dots,  t-1 | \; i_r <i_t   \}$$ 
$$\widetilde{j_t} := j_t - \sharp \{j_r\; \mbox{\em for} \; r=1, \dots,  t-1 | \; j_r <j_t   \}.$$ 
\end{thm}

\begin{proof} By Lemma \ref{lem1}, 
the matrix  $\mu $ is full-rank if and only if
$\tilde{\mu}$ is full-rank and for every 
$i,j$ such that $ \mu_{i,j}= (-\infty, + \infty)$ 
we have that $ Mrk( \mu_{\hat{i}, \hat{j}})
 < p-1$. By Theorem \ref{Mrkquadrate}
this is true if and only if $\tilde{\mu}$ is full-rank and for every 
$i,j$ such that $ \mu_{i,j}= (-\infty, + \infty)$ 
we have that
 in  $\mu_{\hat{i}, \hat{j}} $ 
 there are no totally nonconstant pg-diagonal of  length 
 $p-1$ and the complementary matrix of every   totally nonconstant
 pg-diagonal in  $\mu_{\hat{i}, \hat{j}} $ of  length between $0$ and $p-2$  has $det^c $ equal to $ 0$.

 Moreover, by Lemma \ref{lem2} 
(with $\nu= \tilde{\mu}$), we have that $\tilde{\mu}$ is full-rank
 if and only if, for any $\gamma $ vertex matrix of $\tilde{\mu}$ of even type, we have that $\gamma $  is full-rank and  $det(\overline{\gamma})$ has the same sign as $det( \overline{\tilde{\mu}}_l)$.  
 
 Finally, 
 this is true if and only if (2) holds and, 
 by Lemma \ref{lem3}, for any $\gamma$ vertex matrix of $\tilde{\mu}$ of even type, 
 for any $(i_1,j_1), \ldots, 
 (i_s,j_s) \in \{1,\ldots , p \} \times   \{1,\ldots , p \} $ such that $i_1, \ldots , i_s$ are distinct and $j_1, \ldots , j_s$ are distinct and $\gamma_{i_1,j_1}, \ldots , \gamma_{i_s,j_s} $ are half-bounded, we have that  $det(\overline{\gamma}) \neq 0$
and 
$$ (-1)^{i_1+j_1+ \widetilde{i_2} + \widetilde{j_2} \ldots +\widetilde{i_s}+ \widetilde{j_s}+ \chi((i_1,j_1),
\ldots, (i_s, j_s) ) } \, det (\overline{\gamma})  \,det (\overline{\gamma}_{\widehat{i_1,\ldots, i_s}, \widehat{j_1, \ldots, j_s}} )
  \geq 0, $$
which, by condition (2), is equivalent to condition (3).  
\end{proof}

{\bf Examples.}
1) Let $$\alpha= \begin{pmatrix} 
(-\infty , +\infty ) & [1, + \infty ) & 1 & 1 & 4 \\
1 & [2,3] & 6 & 2 & 4 \\ 
(-\infty ,2] & 0 & [1,4] & 0 & [3,6]\\
0 & [-1, 2] & 3 & 1 & 2 \\
3 & 0  & 3 & 1 & 2
\end{pmatrix}.$$ 
We can easily show that $\alpha$ does not satisfy 
condition (2) of Theorem \ref{gen}, in fact 
$\overline{\tilde{\alpha}}_l =  \begin{pmatrix} 
0 & 1 & 1 & 1 & 4 \\
1 & 2 & 6 & 2 & 4 \\ 
2 & 0 & 1 & 0 & 3\\
0 & -1 & 3 & 1 & 2 \\
3 & 0  & 3 & 1 & 2
\end{pmatrix}$ has negative determinant, while the following vertex matrix of even type of $ \overline{\tilde{\alpha}}$ has positive determinant:  
$$  \begin{pmatrix} 
0 & 1 & 1 & 1 & 4 \\
1 & 2 & 6 & 2 & 4 \\ 
2 & 0 & 1 & 0 & 3\\
0 & 2 & 3 & 1 & 2 \\
3 & 0  & 3 & 1 & 2
\end{pmatrix}$$
So $\alpha$  is not full-rank. In fact it contains the matrix  $$ \begin{pmatrix} 
0 & 1 & 1 & 1 & 4 \\
1 & 2 & 6 & 2 & 4 \\ 
2 & 0 & 1 & 0 & 3\\
0 & 6/5 & 3 & 1 & 2 \\
3 & 0  & 3 & 1 & 2
\end{pmatrix},$$ which is not invertible.

2) Let $$\beta =\begin{pmatrix} [2, +\infty) & 1 & 2 & (-\infty, +\infty) \\ [1,2] & 0 & 3 & 2 \\
3 & [3,7] & 5 & 3 \\ 0 & 0 & 0 & [1, +\infty)
\end{pmatrix}.$$ We can easily see that $\beta $ satisfies conditions (1),(2),(3) of Theorem 
\ref{gen}, so it is full-rank.

3) Let $$\delta =\begin{pmatrix} (-\infty, +\infty) & 1 & 2 & (-\infty, +\infty) \\ [1,2]
&[1,2] & 9 & 2 \\ 3 & [1,5] & 4 & 0 \\
2 & [1,2] & [-1,+ \infty) & 3 \\
\end{pmatrix}.$$ 
Obviously it does not satisfy  condition (1) of Theorem \ref{gen}, in fact  
$\delta_{\hat{1},\hat{1}}$ contains 
totally nonconstant pg-diagonal whose complementary matrix has $det^c \neq 0$.
So $\delta $ is not full-rank.
\bigskip

{\bf Open problem.} A problem that naturally arises is 
the one of the characterization of full-rank matrices whose entries are (not necessarily closed) 
connected subsets of $\R$. 

\bigskip

{\bf Acknowledgments.}
This work was supported by the National Group for Algebraic and Geometric Structures, and their  Applications (GNSAGA-INdAM).

\end{document}